\documentclass{amsart}

\usepackage{amssymb, amsmath, amsthm, hyperref}

\usepackage[english]{babel}

\author{Rafael Torres}

\title[Spin symplectic 4-manifolds with abelian $\pi_1$]{Geography of spin symplectic 4-manifolds with abelian fundamental group}

\address{California Institute of Technology - Mathematics\\ 1200 E California Blvd\\91125\\Pasadena, CA}

\email{rtorresr@caltech.edu}



\theoremstyle{plain}

\newtheorem{theorem}[equation]{Theorem}

\newtheorem{corollary}[equation]{Corollary}

\newtheorem{proposition}[equation]{Proposition}

\newtheorem{lemma}[equation]{Lemma}

\newtheorem{remark}{Remark}

\theoremstyle{definition}

\newcommand{\Z}{\mathbb{Z}}

\begin{document}

\maketitle

In this paper we study the geography and botany of symplectic spin 4-manifolds with abelian fundamental group. Building on the constructions in \cite{[PS]} and \cite{[JP]}, the techniques employed allow us to give alternative proofs and extend their results to the nonsimply connected realm.

\section{Introduction}

Due to the wild and untamed nature of smooth 4-manifolds, efforts towards a classification scheme (existence/uniqueness) take an involved approach. The addition of a symplectic structure has allowed an impressive improvement in our understanding of these objects. The \emph{geography problem}, first proposed by McCarthy and Wolfson in \cite{[MW]}, consists of the existence part of a possible classification: given the Euler characteristic and signature of a spin 4-manifold with a chosen fundamental group, does there exist a symplectic spin  4-manifold with such topological invariants? The (lack of) uniqueness of such manifold is known as the \emph{botany problem}: how many diffeomorphism classes do there exist for the symplectic manifold constructed with the given topological invariants?\\

The geography and botany of irreducible spin simply connected 4-manifolds have been successfully studied in \cite{[FS3], [Sti], [Go], [PS], [JP1]} and \cite{[JP]}, so that most of the existence questions have been settled. The recent addition of Luttinger surgery \cite{[Lu], [ADK]} to the repertoire of symplectic constructions was extremely powerful. Not only did it allow an impressive development in our understanding of simply connected 4-manifolds \cite{[ABBKP], [BK3], [AP]}, but also had as a natural consequence the study of the geography for other fundamental groups \cite{[BK3], [To1], [To2]}. The purpose of this paper is to extend the results on the simply connected realm to manifolds whose fundamental group is an abelian group of small rank.\\

The progress concerning the botany has not been any less sharp. R. Fintushel and R. Stern's work on surgery on nullhomologous tori \cite{[FS1], [FS2]} unveiled a myriad of exotic smooth structures that were previously out of reach through an elegant geometric and topological mechanism. The same authors in joint work with B.D. Park \cite{[FPS]} exploited a duality between Luttinger surgery and its counterpart on nullhomologous tori that enabled the parallel study of the symplectic geography and its botany used by many authors these days, this note included.\\

In order to put the results of this paper into context, we give a rough outline of the current knowledge on the geography of symplectic spin 4-manifolds with $\pi_1 = 1$. In \cite{[PS]}, B.D. Park and S. Szab\'o proved that every allowed homeomorphism type located in the region $0 \leq c_1^2 < 8 \chi_h$ and with odd $b_2^+$ is realized by a simply connected spin irreducible symplectic 4-manifold \cite[Theorem 1.1]{[PS]}. J. Park obtained a similar yet much broader result \cite[Theorem 1.1]{[JP]} which also encompassed spin symplectic simply connected 4-manifolds of zero and positive signature. In particular, he cleverly used a complex spin surface built by C. Persson, C. Peters and G. Xiao in \cite{[PPX]} to produce an infinite number of exotic smooth structures on $(2n + 1)(S^2 \times S^2)$ for a rather large number $n$. Here, $(2n + 1)(S^2\times S^2)$ denotes the connected sum of $2n + 1$ copies of $S^2\times S^2$.\\

Our first result concerns the geography of spin manifolds with negative signature. It provides an extension of B. D. Park and Z. Szab\'o's result to nontrivial abelian fundamental groups. In the simply connected case, we also offer an alternative proof to their theorem.\\

\begin{theorem}{\label{Theorem 1}} Let $s\geq 1$ and let $G$ be either $1, \Z_p, \Z_p \oplus \Z_q$ (and assume $n\geq 2$) or $\Z, \Z \oplus \Z_p, \Z \oplus \Z$ (and $n \geq 1$). For each of the following pairs of integers
\begin{center}
$(c, \chi) = (8n - 8, 2s + n - 1)$,
\end{center}
there exists an irreducible symplectic spin 4-manifold $X$ with 
\begin{center}
$\pi_1(X) = G$ and $(c_1^2(X), \chi_h(X)) = (c, \chi)$.
\end{center}
\end{theorem}

Concerning 4-manifolds with nonnegative signature, by following closely J. Park's main construction in \cite{[JP]} one obtains the following result.

\begin{theorem} {\label{Theorem 2}} Let $G$ be as above. Except for finitely many lattice points, every pair $(c, \chi)$ lying in the region $8\chi \leq c \leq 8.76 \chi$ is realized by an irreducible symplectic spin 4-manifold with 
\begin{center}
$\pi_1(X) = G$ and $(c_1^2(X), \chi_h(X)) = (c, \chi)$.
\end{center}
\end{theorem}

Concerning their botany, we have the following two results.

\begin{proposition} {\label{Proposition 3}} Fix $\pi_1(X) = 1, \Z_p, \Z_q \oplus \Z_q$ or $\Z$, where $q$ is a prime number greater than two. Let $(c, \chi)$ be any pair of integers given in Theorem \ref{Theorem 1} and/or in Theorem \ref{Theorem 2}. There exists an infinite family $\{X_n\}$ of homeomorphic, pairwise nondiffeomorphic irreducible smooth nonsymplectic 4-manifolds realizing the coordinates $(c, \chi)$.
\end{proposition}

\`A la J. Park, for the manifolds with zero signature of Theorem \ref{Theorem 2} we have the following result.

\begin{corollary}{\label{Corollary 4}} There exists an integer $N$ such that $\forall n\geq N$, each of the homeomorphism types given by the manifolds 
\begin{itemize}
\item $(2n + 1) (S^2\times S^2) \# \widetilde{L(p, 1)\times S^1}$,
\item $(2n + 1) (S^2\times S^2) \# \widehat{L(p, 1)\times S^1}$ and
\item $(2n) (S^2\times S^2) \# S^1\times S^3$
\end{itemize}
has infinitely many exotic irreducible smooth structures. In each case, only one of the exotic manifolds is symplectic.\\
\end{corollary}

Here the piece $\widetilde{\L(p, 1)\times S^1}$ stands for the manifold obtained by modifying the product $L(p, 1)\times S^1$  of a Lens space with the circle as follows. Perform a surgery on $L(p, 1)\times S^1$ along $\{x\} \times \alpha$ ($x \in L(p, 1)$) to kill the loop corresponding to the generator of the infinite cyclic group factor so that $\pi_1 = \Z_p$ of the resulting manifold comes from the fundamental group of the Lens space. If instead, we cut out a loop $\{x\} \times \alpha^q$ and glue in a disc to kill the corresponding generator, then we obtain a 4-manifold with $\pi_1 = \Z_p \oplus \Z_q$. Such manifold is denoted by $\widehat{L(p, 1) \times S^1}$.\\

The paper is organized as follows. Section 2 provides the reader with a description of the building blocks and the tools that are employed in our constructions. This section includes the two crucial lemmas for our results as well. In Section 3 we employ them to prove Theorem \ref{Theorem 1} and half of Proposition \ref{Proposition 3}. A description of J. Park' s construction is given in the Fourth section, as well as a proof of Theorem \ref{Theorem 2}, Corollary \ref{Corollary 4} and the remaining part of the proof of Proposition \ref{Proposition 3}.

\section{Tools and Raw materials}

\subsection{Symplectic sums}

In his beautiful paper \cite{[Go]}, R. Gompf introduced \emph{the symplectic sum}, a procedure to build symplectic 4-manifolds which has become essential in our understanding of symplectic 4-manifolds. The following result gathers the properties we will use. 

\begin{lemma}{\label{Lemma 6}} (Gompf, \cite{[Go]}). Let $X$ and $Y$ be spin symplectic 4-manifolds, each containing a symplectic surface $\Sigma_g$ of genus $g$ and self-intersection 0. Then the symplectic sum $X \#_{\Sigma_g} Y$ is a spin symplectic irreducible manifold with coordinates
\begin{center}
$c_1^2(X \#_{\Sigma_g} Y) = c_1^2(X) + c_1^2(Y) + 8(g - 1)$ and\\ $ \chi_h(X\#_{\Sigma_g} Y) = \chi_h(X) + \chi_h(Y) + (g - 1)$.\\
\end{center}
\end{lemma}

The reader is reminded that a spin symplectic 4-manifold is mechanically irreducible, since its Seiberg--Witten invariant is nontrivial \cite{[Ta2], [Ta1]} and it can not be the blow up of another manifold, otherwise it would not be spin.

\subsection{Luttinger surgery and torus surgeries} Carving a torus out of a 4-manifold and then gluing it back in differently is a standard topological procedure to unveil exotic smooth structures. Recently, this idea has been exploited successfully in three directions. First, perform such operation symplectically by adding Luttinger surgery to the palette of constructions of symplectic manifolds; second, use it to construct not only simply connected symplectic manifolds, but also manifolds with several fundamental groups; and last but not least, use a nullhomologous torus that canonically comes out of these surgeries as a dial to change the smooth structure at will. We proceed to give an overview of the machinery. For specific details on the construction, the reader is directed to the references given below.\\

Let $T \subset X$ be a torus of self intersection zero, thus having a tubular neighborhood $N_T \cong T^2 \times D^2$. Let $\alpha$ and $\beta$ be the generators of $\pi_1(T)$ and consider the meridian $\mu_T$ of $T$ in $X$ and the push offs $S^1_{\alpha}, S^1_{\beta}$ in $\partial N_T = T^3$; these are loops homologous in $N_T$ to $\alpha$ and $\beta$ respectively. The manifold obtained from $X$ by performing a $q/p$ - surgery on $T$ along $\beta$ is defined as

\begin{center}
$X_{T, \beta}(q/p) = X - N_T \cup_{\phi} T^2 \times D^2$,
\end{center}

where the gluing map $\phi: T^2\times \partial D^2 \rightarrow \partial (X - N_T)$ satisfies $\phi_*([\partial D^2]) = p[S^1_{\beta}] + q[\mu_T]$ in $H_1(\partial (X - N_T)); \Z)$. Denote core torus $S^1\times S^1 \times \{0\} \subset X_{T, \beta}(q/p)$ by $T_{q/p}$. The surgery reduces $b_1$ by one and $b_2$ by two. The fundamental group of the resulting manifold is given by $\pi_1(X_{T, \beta}(q/p))$.\\

If $X$ is symplectic and $T$ Lagrangian, then performing a $1/p$ surgery on the preferred Lagrangian framing of $N_T$ results in $X_{T, \beta}(1/p)$ being symplectic \cite{[ADK]}. Concerning the botany, in \cite{[FPS]} a procedure is introduced that  uses the nullhomologous torus $T_{q/p}$ to manufacture infinitely many exotic smooth structures starting with a manifold with nontrivial Seiberg--Witten invariant (for example, the symplectic manifold where $T_{q/p}$ was obtained from), by applying a more general $n/1$ - surgery on $T_{q/p}$ (see \cite{[FPS]} or the discussion following \cite[Theorem 13] {[BK3]} for more details). This manufactures an infinite family $\{X_n\}$ of pairwise nondiffeomorphic nonsymplectic 4-manifolds (see Remark \ref{Remark 1}  below).\\

If $X$ is assumed to be spin, one can endow $X_{T, \beta}(q/p)$ with a spin structure by choosing a suitable bundle automorphism $T^2\times D^2 \rightarrow T^2\times D^2$ as follows. Fix a spin structure on $X - N_T$ and one on $T^2\times D^2$. Their difference is given by an element in $H^1(T^2\times D^2; \Z_2) \cong H^1(T^2; \Z_2)$. This element, on the other hand, can be readily seen to be the image of an appropriate bundle automorphism under the coefficient homomorphism $H^1(T^2; \Z) \rightarrow H^1(T^2; \Z_2)$. Thus, identifying two spin structures on $T^2\times D^2$ coming from $X - N_T$ and from $T^2\times D^2$, yields a spin structure for $X_{T, \beta}(q/p)$ itself.\\


We use the remaining part of the section to introduce the building blocks in our constructions.\\

\subsection{Surgeries on $T^4$} This building block will allow us to manipulate the fundamental group of our constructions without changing the Euler characteristic nor the signature. Let $\pi_1(T^4)$ be generated by $x, y, a, b$. Removing a surface from a 4-manifold would normally introduce more generators to the fundamental group of the complement. In \cite{[BK3]}, S. Baldridge and P. Kirk showed that the fundamental group of the complement of two Lagrangian tori $T_1$ and $T_2$ inside the 4-torus is generated by 4 elements, just like $\pi_1(T^4)$ itself.\\

\begin{proposition}{\label{Proposition 7}}(Baldridge--Kirk, \cite{[BK3]}) The fundamental group of $T^4 - (T_1 \cup T_2)$ is generated by the loops $x, y, a, b$ and the relations $[x, a] = [y, a] = 1$ hold. The meridians of the tori and the two Lagrangian push offs of their generators are given by the following formulae:\\
\begin{center}
$\mu_1 = [b^{-1}, y^{-1}], m_1 = x, l_1 = a$ and\\
$\mu_2 = [x^{-1}, b], m_2 = y, l_2 = bab^{-1}$.\\
\end{center}
\end{proposition}

As a corollary of their efforts one obtains the following lemma.

\begin{lemma}{\label{Lemma 8}} Let $X$ be a simply connected spin symplectic 4-manifold containing a symplectic torus such that $\pi_1(X - T) = 1$. There exists a spin symplectic 4-manifold with Chern numbers $\chi_h(Z) = \chi_h(X)$ and $c_1^2(Z) = c_1^2(X)$. The fundamental group of $Z$ can be chosen to be
\begin{enumerate}
\item $\pi_1 = \Z \oplus \Z$,
\item $\pi_1 = \Z \oplus \Z_q$,
\item $\pi_1 = \Z$ 
\end{enumerate}
\end{lemma}

\begin{proof} Let $T_1 \subset T^4$ be as above. Perturb the symplectic form on $T^4$ such that $T_1$ becomes symplectic while $T_2$ stays Lagrangian (see \cite[Lemma 1.6]{[Go]}). The torus $T_1$ carries the generators $x$ and $b$. Take the symplectic sum $Y: = T^4 \#_{T_1 = T} X$. Since the meridian of $T$ in $X - T$ is trivial, the relation $[y, b] = 1$ holds in the fundamental group of this newly constructed manifold. Therefore, the symplectic sum results in a manifold $Y$ with $\pi_1(Y) = \Z y \oplus \Z b$. We can now proceed to apply a $1/q$ Luttinger surgery to $T_2$ to produce a manifold with $\pi = \Z_p \oplus \Z b$; for $q = 1$ we have $\pi_1 = \Z$ and for $q > 1$, $\pi_1 = \Z_q \oplus \Z b$.\\
\end{proof}

\subsection{Manifolds with the cohomology of $(2n - 3) (S^2\times S^2)$}

In \cite{[FPS]}, R. Fintushel, B.D. Park and R. Stern built an infinite family of irreducible pairwise nondiffeomorphic spin 4-manifolds with the same integer cohomology ring as $S^2 \times S^2$. A. Akhmedov and B.D. Park generalized the construction in \cite{[AP]}, by producing an infinite family of irreducible pairwise nondiffeomorphic spin 4-manifolds $\{Y_n(m)| m = 1, 2, 3, \ldots \}$ with only one symplectic member which has the same integer cohomology ring as $(2n - 3)(S^2\times S^2)$ with $n\geq 2$. The characteristic numbers of these manifolds are $e = 4n - 4$ and $\sigma = 0$; equivalent, $\chi_h = n - 1$ and $c_1^2 = 8n - 8$.\\ 


These manifolds are constructed by applying $2n + 3$ Luttinger surgeries and one torus surgery to $\Sigma_2 \times \Sigma_n$ (the product of a genus 2 surface with a genus n surface). Let $a_i, b_i, c_j$ and $d_j$ ($i = 1, 2$, $j = 1, \ldots, n$) be the standard generators of $\pi_1(\Sigma_2)$ and $\pi_1(\Sigma_n)$ respectively. The following relations hold in $\pi_1(Y_n(m))$. We refer the reader to \cite{[AP]} for further details.\\

\begin{center}
$[b_1^{-1}, d_1^{-1}] = a_1, [a_1^{-1}, d_1] = b_1, [b_2^{-1}, d_2^{-1}] = a_2, [a_2^{-1}, d_2] = b_2$,\\
$[d_1^{-1}, b_2^{-1}] = c_1, [c_1^{-1}, b_2] = d_1, [d_2^{-1}, b_1^{-1}] = c_2, [c_2^{-1}, b_1] = d_2$,\\
$[a_1, c_1] = 1, [a_1, c_2] = 1, [a_1, d_2] = 1, [b_1, c_1] = 1$,\\
$[a_2, c_1] = 1, [a_2, c_2] = 1, [a_2, d_1] = 1, [b_2, c_2] = 1$,\\
$[a_1, b_1][a_2, b_2] = 1$, $[c_1, d_1][c_2, d_2] = 1$.\\
\end{center}

and\\

\begin{center}
$[a_1^{-1}, d_3^{-1}] = c_3, [a_2^{-1}, c_3^{-1}] = d_3, \cdots, [a_1^{-1}, d_n^{-1}] = c_n, [a_2^{-1}, c_n^{-1}] = d_n$,\\
$[b_1, c_3] = 1, [b_2, d_3] = 1, \cdots, [b_1, c_n] = 1, [b_2, d_n] = 1$,\\
$\prod_{j = 2}^n[c_j, d_j] = 1$.\\
\end{center}

These manifolds are our basic building block for manipulating the fundamental group. We employ them to obtain the following result.\\
 
\begin{lemma}{\label{Lemma 9}} Let $X$ be a simply connected spin symplectic 4-manifold containing a symplectic torus such that $\pi_1(X - T) = 1$. Then for all $n\geq 1$ there exists a spin symplectic 4-manifold with Chern numbers $\chi_h(Z) = \chi_h(X) + n - 1$ and $c_1^2(Z) = c_1^2(X) + 8n - 8$. The fundamental group of $Z$ can be chosen to be
\begin{enumerate}
\item $\pi_1 = \Z \oplus \Z$,
\item $\pi_1 = \Z \oplus \Z_q$,
\item $\pi_1 = \Z_p \oplus \Z_q$,
\item $\pi_1 = \Z_q$,
\item $\pi_1 = \Z$ or
\item $\pi_1 = 1$.
\end{enumerate}
Furthermore, $Z$ contains a Lagrangian torus such that the inclusion induced homomorphism $\pi_1(Z - T) \rightarrow \pi_1(Z)$ is an isomorphism.
\end{lemma}

\begin{proof} Consider the case $n = 2$. Let $S$ be the manifold obtained by applying 5 Luttinger $\pm 1$ -surgeries to $\Sigma_2\times \Sigma_2$. The surgeries that are not to be performed are $(a_1'\times c_1', a_1', -1)$, $(a_2'\times c_2', a_2', -1)$ and $(a_2''\times d_1', d_1', +1)$. Call these three tori $T_1, T_2$ and $T_3$ respectively. In $\pi_1(S)$ all the relations from $\pi_1(Y_2(1))$ hold except for $[b_1^{-1}, d_1^{-1}] = a_1$, $[b_2^{-1}, d_2^{-1}] = a_2$ and $[c_2^{-1}, b_1] = d_2$.\\

Build the symplectic sum of $X$ and $S$ along the corresponding torus in $X$ and $T_1$ in $S$  and call the resulting manifold $S_{\Z \oplus \Z}$. The meridian of $T_1$, $[b_1^{-1}, d_1^{-1}] = a_1$ is killed during the symplectic sum and the surviving relations show that $\pi_1(S_{\Z \oplus \Z} - T_2 \cup T_3)$ is generated by the two commuting elements $a_2$ and $d_1$. The Mayer--Vietoris sequence shows that $H_1(S_{\Z \oplus \Z} - T_2 \cup T_3); \Z) = \Z^2$, thus $\pi_1(S_{\Z \oplus \Z}) = \Z a_2 \oplus \Z d_1$. It is straight forward to check $e(S_{\Z \oplus \Z}) = e(X) + 4$ and $\sigma(S_{\Z \oplus \Z}) = \sigma(X)$.\\

Notice that the geometrically dual torus $T'$ to $T_1$ is contained in $S_{\Z \oplus \Z}$ and its meridian is trivial in the complement. This implies $\pi_1(S_{\Z \oplus \Z} - T') \cong \pi_1(S_{\Z \oplus \Z}) = \Z^2$. Thus, item (1) of the lemma has been produced.\\

Applying $(a_2'\times c_2', a_2', -1/q)$, that is a $-1/q$ Luttinger surgery to $S_{\Z\oplus \Z}$ on $T_2$ along $a_2'$, produces item (2). By applying $(a_2''\times d_1', d_1', +1/p)$ to the resulting manifold one produces item (3) ($p > 1$) and item (4) ($p = 1$). Applying $(a_2''\times d_1', d_1', +1)$ to $S_{\Z \oplus \Z}$ produces item (5), while item (6) on the list is produced by applying both surgeries $(a_2''\times d_1', d_1', +1)$ and $(a_2'\times c_2', a_2', -1)$ to $S_{\Z \oplus \Z}$.\\

The cases $n \geq 3$ follow almost verbatim to the procedure described above substituting $\Sigma_2 \times \Sigma_2$ with $\Sigma_2 \times \Sigma_n$. The details are left to the reader. We do point out that the bigger $n$ is, the more Lagrangian tori the resulting manifold contains. For example, the manifold obtained by applying Luttinger surgeries to $\Sigma_2 \times \Sigma_5$ contains 12 Lagrangian tori while the one obtained from $\Sigma_2\times \Sigma_7$ has 20 Lagrangian tori; all of them have trivial meridian.
\end{proof}

\begin{remark}{ \label{Remark 1}} An infinite number of exotic smooth structures can be unveiled at the prize of surrendering the symplectic structure. We exemplify the process in the infinite cyclic fundamental group case, while the next paragraph explains why the process works for the rest of the groups. Before applying the last Luttinger surgery to obtain a symplectic manifold with $\pi_1 = \Z$, one has a symplectic manifold $X_{\Z\oplus \Z}$ with $\pi_1(X_{\Z \oplus \Z}) = \Z \oplus \Z$. By Taubes' results \cite{[Ta1], [Ta2]}, $X_{\Z\oplus \Z}$ has nontrivial Seiberg--Witten invariants. One can performed a more general torus surgery on $X_{\Z \oplus \Z}$, to obtain a manifold $X_{\Z}$, with infinite cyclic fundamental group, and nontrivial Seiberg--Witten invariants. The manifold $X_{\Z}$ contains a nullhomologous torus $T'$. Applying a torus surgery on $T'$ yields an infinite family $\{X_n\}$ parametrized by the surgery coefficient $n$. The formula given in \cite[Theorem 3.4]{[MMS]} can be used to prove that the Seiberg--Witten invariants distinguish infinitely many diffeomorphism types within the members of $\{X_n\}$ (see also \cite[Corollary 2]{[FPS]}).\\ 

To conclude on their homeomorphism type, one must check that these manifolds have the desired fundamental group; we already know their Chern invariants remained unchanged after the surgery. For this purpose it suffices to see that the effect such surgery has on the presentation of the fundamental groups is to replace a relation of the form $[a, b] = c^p$ by $[a, b]^n = c^p$ for a given $p$ and $n$ and generators $a, b$. Given that in the proofs of Lemma \ref{Lemma 8} and Lemma \ref{Lemma 9} we concluded that the original relation is trivial, then raising it to any power will result in a trivial relation as well. Hence, we will make no distinctions in future sections about the computations of $\pi_1$ of the infinite families.
\end{remark}

\subsection{Horikawa surfaces} The complex surfaces satisfying $c_1^2 = 2\chi_h - 6$ are commonly known as Horikawa surfaces and are denoted by $H(4k - 1)$. They are constructed as branched covers of the Hirzebruch surface $\mathbb{F}_{2m}$ along disconnected curves and we point out that a simply connected Horikawa surface is spin if and only if $k$ is even. The Chern invariants of the specific manifolds we will be using, $H(8k' - 1)$, are given by $(c_1^2, \chi_h) = (16 k' - 8, 8k' - 1)$. Moreover, $H(8k' - 1)$ contains an embedded Lagrangian torus which intersects a 2-sphere transversally at one point \cite{[FS3], [Sti]}.\\


\subsection{A spin surface of positive signature}

In \cite{[PPX]}, U. Persson, C. Peters and G. Xiao constructed a simply connected spin complex surface $Y$ of positive signature which contains a holomorphic curve $\Sigma_g$ of genus g and trivial self intersection. Furthermore, the meridian of this surface in the complement is trivial since $Y$ also contains an embedded 2-sphere $\mathbb{CP}^1$ intersecting $\Sigma_g$ transversely at a point. Its Chern invariants are approximately $\chi_h(Y) \approx 6857 x^2$ and $c_1^2(Y) \approx 60068 x^2$. 

\subsection{Knot surgery on elliptic minimal surfaces} Our last building block is also a classical element in the construction of 4-manifolds and we only remind the reader of its properties relevant to our purposes. Let $E(2s)$ denote the underlying smooth 4-manifold of the simply connected minimal elliptic surface without multiple fibers and with geometric genus $p_g = 2s - 1$ \cite{[G1]}  and \cite[Prop. 3.1.11] {[GSt]}. Its Chern numbers are given by $c_1^2 = 0$ and $\chi_h = 2s$. Notice that in particular $E(2)$ is a $K3$ surface. In Section 3 and Section 4, it is easy to see where the manifold $E(2s)$ can be replaced by an exotic version $E(2s)_K$ obtained by Knot surgery\cite{[FS1]}.\\


\section{Negative signature}{\label{Section 3}}

\subsection{Examples with $\sigma =-16 s$ for $s > 0$}

\begin{proposition}{\label{Proposition 10}} Let $s\geq 1$. For $\pi = 1, \Z_p, \Z_p \oplus \Z_q$ assume $n\geq 2$ and for $\pi = \Z, \Z \oplus \Z_p$ and $\Z \oplus \Z$ assume $n\geq 1$. There exists a spin irreducible symplectic manifold $X$ satisfying $c_1^2 = 8 n - 8$, $\chi_h = n + 2s - 1$ and $\pi_1(X) = \pi$.
\end{proposition}

\begin{proof} The proposition follows from employing $X = E(2s)_K$ in Lemma \ref{Lemma 8} and Lemma \ref{Lemma 9}.
\end{proof}



By applying the corresponding homeomorphism criteria, we conclude that the manifolds constructed are homeomorphic to the following topological prototypes: 
\begin{itemize}
\item $\pi = 1$: $E(2s)\# (2n - 2) (S^2\times S^2)$.\\
\item $\pi = \Z_p$: $E(2s) \# (2n - 2)(S^2\times S^2) \# \widetilde{L(p, 1)\times S^1}$.\\
\item $\pi = \Z_q \oplus \Z_q$: $E(2s) \# (2n - 2)(S^2\times S^2) \# \widehat{L(p, 1)\times S^1}$.\\
\item $\pi = \Z$: $E(2s) \# (2n - 1)(S^2\times S^2) \# S^3\times S^1$.\\
\end{itemize}

Indeed, notice that the Euler characteristic, the spin property and the signature of the symplectic sum are computed by Lemma \ref{Lemma 6}. Moreover, torus surgeries do not change any of these topological invariants. In the simply connected case, the known result of Freedman's \cite{[MF]} establishes the homeomorphism type of the manifolds constructed. Less known, yet outstanding results allow us to conclude on the homeomorphism types for manifolds with nontrivial fundamental group by using the same topological invariants. The criteria of Hambleton--Teichner in \cite{[HT]} concludes on the homeomorphism type for manifolds with infinite cyclic fundamental group. The finite fundamental group cases, both cyclic and noncyclic yet abelian of odd order, follow from the criteria of Hambleton--Kreck in \cite{[HK]} by checking that the manifolds constructed share the same $\omega_2$-type: this requires to see that the universal cover of the manifolds are spin as well. Notice the need of the hypothesis requiring $q$ to be a prime number in Proposition \ref{Proposition 3}.\\

Thus, considering Remark 1 we have the following result.

\begin{corollary} Each of the manifolds
\begin{itemize}
\item$E(2s)\# (2n - 2)(S^2\times S^2)$,
\item$E(2s) \# (2n - 2)(S^2\times S^2) \# \widetilde{L(p, 1)\times S^1}$,
\item $E(2s) \# (2n - 2)(S^2\times S^2) \# \widehat{L(p, 1)\times S^1}$ and
\item $E(2s) \# (2n - 1)(S^2\times S^2) \# S^3\times S^1$.
\end{itemize}
admits infinitely many exotic irreducible smooth structures. In each case, only one of these exotic manifolds is symplectic.\\
\end{corollary}

These methods improve the main theorem in \cite{[FS]}, where R. Fintushel and R. Stern constructed a manifold $X$ homeomorphic to $K3 \# S^2\times S^2 \# S^3\times S^1$. We remind the reader that in the abelian, yet noncyclic case, the fundamental group is assumed to be $\pi_1 = \Z_q \oplus \Z_q$, where $q$ is a prime number.

\subsection{Examples with $\sigma = -48k'$ for $k' > 0$}

Employing the Horikawa surfaces $H(8k' - 1)$ and $H(7) \#_{T = T} \# H(8k' - 1)$ in Lemma \ref{Lemma 8} and Lemma \ref{Lemma 9} yields the following proposition. 

\begin{proposition} {\label{Proposition 12}}Let $k' > 0$. For $\pi = 1, \Z_p, \Z_p \oplus \Z_q$ assume $n\geq 2$ and for $\pi = \Z, \Z \oplus \Z_p$ and $\Z \oplus \Z$ assume $n\geq 1$. There exists a spin irreducible symplectic manifold $X$ satisfying 
\begin{itemize}
\item $c_1^2(X) = 16k' + 8n - 16, \chi_h(X) = 8k' + n - 2$ or
\item $c_1^2(X) = 16k' + 8n + 88, \chi_h(X) = 8k' + n + 53$
\end{itemize}
and $\pi_1(X) = \pi$.
\end{proposition}

\begin{corollary} Each of the manifolds
\begin{itemize}
\item$H(8k' - 1) \# (2n - 2)(S^2\times S^2)$,
\item$H(8k' - 1) \# (2n - 2)(S^2\times S^2) \# \widetilde{L(p, 1)\times S^1}$,
\item $H(8k' - 1) \# (2n - 2)(S^2\times S^2) \# \widehat{L(p, 1)\times S^1}$ and
\item $H(8k' - 1) \# (2n - 1)(S^2\times S^2) \# S^3\times S^1$.
\end{itemize}
admits infinitely many exotic irreducible smooth structures. In each case, only one of these exotic manifolds is symplectic.
\end{corollary}



\section{Nonnegative signature}

\subsection{J. Park's construction} In \cite{[JP]}, J. Park used the spin complex surface described in 2.6 above to realize all but finitely many allowed points in the region $0 \leq c_1^2 \leq 8.74 \chi_h$ for trivial fundamental group. Given that we already filled in the points of negative signature above, we now follow his construction in \cite{[JP]} almost verbatim in order to address the region $8 \leq c_1^2 \leq 8.76 \chi_h$. We start by describing the argument and main building blocks in \cite{[JP]}.\\

Consider a simply connected spin symplectic 4-manifold $Z$ which contains a symplectic torus $T$ in a cusp neighborhood $N$ and symplectic surface $\Sigma_g$ of genus g and zero self intersection, $\Sigma_g$ disjoint from $N$. The Chern invariants of this manifold are $c_1^2(Z) = 8g^2 - 16 g + 8$ and $\chi_h(Z) = 2 g^2 - g + 1$. In particular its signature is given by $\sigma(Z) = -8g^2 + 8g$. Now take the spin complex surface described in Section 2.6 and build the symplectic sum

\begin{center}
$X: = \overbrace{Y \#_{\Sigma_g} \cdots \#_{\Sigma_g} Y}^{k} \#_{\Sigma_g} Z$.\\
\end{center}

Assume the integer $k$ is such that $X$ has positive signature. Furthermore, $\pi_1(X) = 1$ since all the pieces are simply connected and the meridian of $\Sigma_g$ in $Y - \Sigma_g$ is trivial. The Chern numbers can be calculated to be $c_1^2(X) = k c_1^2(Y) + c(Z) + 8k(g - 1)$ and $\chi_h(Y) = k \chi_h(Y) + \chi_h(Z) + k(g - 1)$, thus by considering large enough integers $k$ and $x$, one has\\

\begin{center}
$\frac{c_1^2(X)}{\chi_h(X)} = \frac{k c_1^2(Y) + c(Z) + 8k(g - 1)}{k \chi_h(Y) + \chi_h(Z) + k(g - 1)} \approx \frac{c_1^2(Y)}{\chi_h(Y)} \approx \frac{60068 x^2}{6857 x^2} = 8.76009 \cdots$\\
\end{center}

J. Park then fixes $k$ and $x$ big enough such that $\frac{c_1^2(X)}{\chi_h(X)} > 8.76$ holds. At this point one should notice that $X$ contains a symplectic torus of self intersection zero lying on the building block $Z$. In fact, one can also find such tori in the $Y$ blocks. To finish his argument, he then proceeds to define a line $c = f(\chi)$

\begin{center}
$f(\chi) = \frac{c(X)}{\chi(X)} \cdot (\chi - c(X) / 2 - 6) + c(X)$
\end{center}

whose slope $\frac{c(X)}{\chi(X)} = \frac{c_1^2(X)}{\chi_h(X)}$ is greater than 8.76. Finally, build the simply connected manifold $W:= \overbrace{X\#_{T^2} X\#_{T^2} \#_{T^2} \cdots \#_{T^2} X}^{m} \#_{T^2} V$ (where the block $V$ can be chosen from $H(8k' - 1)\#_{T^2} E(2s), H(7)\#_{T^2}H(8k' - 1)\#_{T^2} E(2s)$ or a simply connected manifold constructed in Proposition \ref{Proposition 10}). Then, for some integer $m$, for every allowed lattice point $(c, \chi)$ in the first quadrant of the geography plane which complies with $c = f(\chi)$, there exists an irreducible symplectic simply connected spin 4-manifold $W$ with $(c, \chi) = (c_1^2(W), \chi_h(W))$.\\
 
Given that $W$ has a torus $T$ of self intersection zero and of trivial meridian in $W - T$, Lemma \ref{Lemma 8} and Lemma \ref{Lemma 9} imply the following result (Theorem \ref{Theorem 2}).

\begin{proposition} Let $\pi = 1, \Z_p, \Z_p \oplus \Z_q, \Z, \Z \oplus \Z_p$ and $\Z \oplus \Z$. Except for finitely many lattice points, for every allowed pair $(c, \chi)$ lying in the region
\begin{center}
$8\chi \leq c \leq 8.76 \chi$,
\end{center}
there exists a spin irreducible symplectic manifold $X$ satisfying 
\begin{center}
$\pi_1(X) = \pi$ and $(c_1^2(X), \chi_h(X) = (c, \chi)$.
\end{center}
\end{proposition}

Concerning the manifolds with negative signature from the previous proposition, we have the following result.


\begin{corollary} There exists an integer $N$ such that $\forall$  $n\geq N$, each of the manifolds
\begin{itemize}
\item $(2n + 1) (S^2\times S^2)$,
\item $(2n + 1) (S^2\times S^2) \# \widetilde{L(p, 1)\times S^1}$,
\item $(2n + 1) (S^2\times S^2) \# \widehat{L(p, 1)\times S^1}$ and
\item $(2n + 2) (S^2\times S^2) \# S^1\times S^3$
\end{itemize}
has infinitely many exotic irreducible smooth structures. For each case, only one of the exotic manifolds admits a symplectic structure.
\end{corollary}


\end{document}